\documentclass{amsart}

\usepackage[utf8]{inputenc}
\usepackage[T1]{fontenc}

\usepackage{amsmath,amsfonts,amssymb,amsthm}
\usepackage{mathtools,extpfeil,xfrac}
\usepackage{csquotes}
\usepackage{microtype}
\usepackage{tikz-cd}

\RequirePackage{enumitem}
\setlist[enumerate,1]{label=(\roman*)}
\setlist[enumerate,2]{label=(\alph*)}
\setlist[enumerate,3]{label=(\Alph*)}
\setlist[enumerate,4]{label=\arabic*.}

\RequirePackage[
    backend=biber,
    style=alphabetic,
    maxnames=100,
    maxalphanames=100,
]{biblatex}

\DeclareLabelalphaTemplate{
    \labelelement{
        \field[final,lowercase]{shorthand}
        \field{label}
        \field[strwidth=1,strside=left]{labelname}
    }
}

\DeclareFieldFormat{extraalpha}{#1}

\usepackage[colorlinks=false,unicode=true]{hyperref}
\usepackage{bookmark}

\newtheorem{theorem}{Theorem}[section]
\newtheorem{proposition}[theorem]{Proposition}
\newtheorem{lemma}[theorem]{Lemma}
\newtheorem{corollary}[theorem]{Corollary}
\theoremstyle{definition}
\newtheorem{definition}[theorem]{Definition}
\theoremstyle{remark}
\newtheorem{remark}[theorem]{Remark}
\newtheorem{example}[theorem]{Example}

\let\shortbar\bar
\let\bar\overline

\newcommand\as[2][]{\mathbb A^{#2}_{#1}}    
\newcommand\ps[2][]{\mathbb P^{#2}_{#1}}    
\newcommand\pt{\mathrm{pt}}                 
\newcommand\Ga{\mathbb{G}_{\mathrm{a}}}     
\newcommand\Gm{\mathbb{G}_{\mathrm{m}}}     
\DeclareMathOperator\Spec{Spec}

\let\stack\mathbf                           
\newcommand\cs{\stack{B}}                   
\newcommand\LocSys{\operatorname{LocSys}}   
\newcommand\Bun{\operatorname{Bun}}         

\newcommand\cat{\mathbf}                    
\newcommand\catCoSimplicial{\Delta^{\mathrm{op}}}
\DeclareMathOperator\Hom{Hom}
\newcommand\id[1][]{{\mathrm{Id}_{#1}}}     
\newcommand\catVect[1][]{\cat{Vect}_{#1}}   
\newcommand\catModules[1]{\cat{Mod}(#1)}    
\newcommand\cx\bullet                       
\newcommand\dual{\vee }
\newcommand\cocolon{
    \nobreak \mskip 6muplus1mu \mathpunct {}\nonscript \mkern -\thinmuskip{:}\mskip 2mu\relax%
}
\newcommand\rquot[2]{
    \mathchoice%
        {\left.#1\kern-0.2ex\middle/\kern-0.3ex\lower0.7ex\hbox{$\displaystyle #2$}\right.}%
        {\left.#1\middle/#2\right.}%
        {\left.#1\middle/#2\right.}%
        {\left.#1\middle/#2\right.}%
}
\newcommand\isoto{\xrightarrow{\sim}}       
\newcommand\proj[1]{\operatorname{pr}_{#1}} 
\newcommand{\HCoh}{\operatorname{HH}^\cx}   

\newcommand\sheaf\mathcal
\DeclareMathOperator\supp{supp}     
\newcommand\catDMod[2][]{\cat{DMod}_{#1}(#2)}   
\newcommand\catDModHol[1]{\catDMod[\mathrm{hol}]{#1}}   
\newcommand\dR{\mathrm{dR}}
\newcommand\GammadR{\Gamma _{\mkern-4mu\dR}}

\newcommand\dualize{\mathbb D}
\newcommand\catIndCoh[2][]{\cat{IndCoh}_{#1}(#2)} 

\newcommand\opalg[1]{#1^{\mathrm{op}}}
\newcommand\Gammasub[1]{\Gamma_{\mkern-3mu#1}}
\newcommand\barGammasub[1]{\bar{\Gamma}_{\mkern-3mu#1}}

\newcommand\XXtimes{%
    \mathchoice%
        {\mathop{\times }\limits_{\stack X \times  \stack X}}
        {\times _{\stack X \times  \stack X}}
        {\times _{\stack X \times  \stack X}}
        {\times _{\stack X \times  \stack X}}
    }

\newcommand\bc{\textbf{($\mathbf{*}$)}}
\newcommand\hbc{\textbf{($\mathbf{**}$)}}

\newcommand\catK[2][]{\cat{K}_{#1}(#2)}
\newcommand\catKHol[1]{\catK[\mathrm{hol}]{#1}}

\newcommand\ls[1]{\mathcal{L} #1}
\newcommand\lsY[2][\stack Y]{\mathcal{L}_{#1} #2}
\newcommand\cls[1]{\overline{\mathcal{L}} #1}
\newcommand\clsY[2][\stack Y]{\overline{\mathcal{L}}_{#1} #2}
\newcommand\lsc[1]{\mathcal{L}^c #1}
\newcommand\lscY[2][\stack Y]{\mathcal{L}_{#1}^c #2}

\newcommand\sls[1]{L#1}
\newcommand\slsY[2][\stack Y]{L_{#1}#2}
\newcommand\scls[1]{\overline{L}#1}
\newcommand\sclsY[2][\stack Y]{\overline{L}_{#1}#2}
\newcommand\slsc[1]{L^c #1}
\newcommand\slscY[2][\stack Y]{L^c_{#1}#2}

\newcommand\goodness{property \textbf{(!)}}

\newcommand\isgood{has property \textbf{(!)}}
\newcommand\aregood{have property \textbf{(!)}}

\title{Hochschild cohomology of torus equivariant D-modules}

\author{Clemens Koppensteiner}

\begin{document}

\begin{abstract}%
  We discuss the Hochschild cohomology of the category of D-modules associated to an algebraic stack.
  In particular we describe the Hochschild cohomology of the category of torus-equivariant D-modules as the cohomology of a D-module on the loop space of the quotient stack.
  Finally, we give an approach for understanding the Hochschild cohomology of D-modules on general stacks via a relative compactification of the diagonal.
  This work is motivated by a desire to understand the support theory (in the sense of \cite{BensonIyengarKrause:2008:LocalCohomologyAndSupportForTriangulatedCategories}) of D-modules on stacks.
\end{abstract}

\maketitle

\setcounter{tocdepth}{1}
\tableofcontents

\section{Introduction}

Given a manifold $X$ and a category of sheaves on $X$, microlocal geometry asks whether the sheaves can be localized not just on $X$ but also with respect to codirections, i.e.~on the cotangent space $T^*X$.
For example, for constructible sheaves this leads to the notion of microsupport discussed in detail in \cite{KashiwaraSchapira:1994:SheavesOnManifolds}.
More generally, given a category of sheaves on a space $X$, we can ask whether it is possible to localize them on some space that is strictly larger than $X$ itself.

Even more generally one can ask the following question: Given a $k$-linear category $\cat C$, can one find a space over which $\cat C$ localizes?
For co-complete compactly generated triangulated categories one answer is provided by Benson, Iyengar and Krause \cite{BensonIyengarKrause:2008:LocalCohomologyAndSupportForTriangulatedCategories}:
To each map from a graded-commutative ring $R$ to the center of $\cat C$ one associates the \emph{triangulated support} functor $\supp_R$, assigning to each object $A \in  \cat C$ a subset $\supp_R A \subseteq \Spec R$.
This construction can be used to unify various theories of support in different areas of mathematics (though it does not yield the microlocal support of constructible sheaves mentioned above).

We are led to consider the universal algebra acting on $\cat C$ by this construction, i.e.~the \emph{Hochschild cohomology} of $\cat C$.
For a complete pre-triangulated dg category $\cat C$ the Hochschild cohomology is the dg algebra of derived endomorphisms of the identity functor of $\cat C$:
\[
    \HCoh(\cat C)
    = \mathbb{R}\mathrm{Hom}(\id[\cat C], \id[\cat C])
    = \Hom_{\cat{Funct}(\cat C, \cat C)}(\id[\cat C], \id[\cat C]).
\]
The ring $R = \bigoplus \operatorname{HH}^{2n}(\cat C)$ is commutative and hence one can define for each $A \in  \cat C$ the support $\supp_R A$ as a subset of $\Spec R$.
Thus understanding the Hochschild cohomology of a dg category can be an important step to understanding the category itself.

This construction, applied to the category of (ind-)coherent sheaves on a (quasi-smooth, dg-) scheme, yields the singular support of coherent sheaves introduced by Arinkin and Gaitsgory \cite{ArinkinGaitsgory:2015:SingularSupport}.
The authors use this notion of singular support for the category $\catIndCoh{\LocSys_G}$ in their formulation of the geometric Langlands conjecture.
By Langlands duality, one should then have a matching support theory for the category $\catDMod{\Bun_G}$ and the question arises whether it is possible to formulate this theory in a way that is intrinsic to D-modules.

A first step to this -- and also a problem of independent interest -- is to understand the Hochschild cohomology of the category $\catDMod{\stack X}$ of D-modules on a stack $\stack X$.
Thanks to the theory of integral kernels, this question can be approached geometrically.
Representing the identity functor by its kernel, we obtain an isomorphism of dg algebras
\begin{equation}
    \label{eq:intro:d-mod-hcoh}
    \HCoh(\catDMod{\stack X}) \cong \Hom_{\catDMod{X\times X}}(\Delta _*\omega _{\stack X}, \Delta _*\omega _{\stack X}),
\end{equation}
where $\Delta \colon \stack X \to \stack X \times \stack X$ is the diagonal morphism and $\omega _{\stack X}$ is the dualizing module.
In particular if $\stack X$ is a separated scheme, then $\Delta $ is a closed embedding and $(\Delta ^*,\Delta _*)$ adjunction combined with Kashiwara's Equivalence show that $\opalg{\HCoh(\catDMod{\stack X})}$ is isomorphic to the de~Rham cohomology of $\stack X$.
However, if $\stack X$ is a general stack, then $\Delta $ is not proper and the situation becomes more complicated.

By Verdier duality and adjunction we can always rewrite \eqref{eq:intro:d-mod-hcoh} as
\[
    \HCoh(\catDMod{\stack X}) \cong
    \opalg{\Hom_{\catDMod{X}}(k_{\stack X}, \Delta ^!\Delta _!k_{\stack X})} =
    \opalg{\GammadR(\stack X, \Delta ^!\Delta _!k_{\stack X})}.
\]
It is now tempting to look at the Cartesian square
\[
    \begin{tikzcd}
        \ls \stack X \arrow{r}{p_2} \arrow{d}{p_1} & \stack X \arrow{d}{\Delta } \\
        \stack X \arrow{r}{\Delta } & \stack X \times  \stack X
    \end{tikzcd}
\]
where
\[
    \ls \stack X = \stack X \mathop{\times }\limits_{\stack X \times  \stack X} \stack X
\]
is the derived loop space (or inertia stack) of $\stack X$ and try to express the Hochschild cohomology as the cohomology of some sheaf on $\ls \stack X$.
We could expect the existence of an isomorphism
\begin{equation}
    \label{eq:intro:naive}
    \GammadR(\stack X, \Delta ^!\Delta _!k_{\stack X}) \cong
    \GammadR(\stack X, p_{1,!}p_2^! k_{\stack X}).
\end{equation}
We refer to this as the \enquote{naive expectation}.
Unfortunately, the two sides are in general not isomorphic (the stack $\stack X = \ps1/\Ga$ is an easy counter-example).

In this paper we investigate the morphism
\begin{equation}\label{eq:intro:central_cone}
    p_{1,!}p_2^! k_{\stack X} \to  \Delta ^!\Delta _! k_{\stack X}
\end{equation}
and thus whether the naive isomorphism \eqref{eq:intro:naive} holds.
As an application, we prove the following theorem, giving a class of stacks where \eqref{eq:intro:naive} is indeed an isomorphism.

\begin{theorem}\label{thm:main}
    Let $G \cong \Gm^n$ be a torus acting on a normal quasi-projective variety $X$ over an algebraically closed field $k$ of characteristic $0$.
    Then there is a canonical isomorphism of algebras
    \[
        \HCoh\bigl(\catDMod{X/G}\bigr)
        \cong
        \opalg{\GammadR\bigl(X/G,\,p_{1,!} p_2^! k_{X/G}\bigr)},
    \]
    where the algebra structure on $\GammadR\bigl(X/G,\,p_{1,!} p_2^! k_{X/G}\bigr)$ is induced by the groupoid structure on $\ls(X/G)$.
\end{theorem}

\subsection*{Contents}

Since, to the author's knowledge, there is currently no comprehensive account of the theory of D-modules on stacks available in the literature, we review the basic formalism in Section~\ref{sec:d-mods}.
In essence, most of the familiar six functor formalism of holonomic D-modules works in a general context.
It is convenient to rephrase Theorem~\ref{thm:main} as an equivalence of monads.
Thus in Section~\ref{sec:monads} we review some constructions of monads in our setting.
In Section~\ref{sec:HH} we explain how morphism~\eqref{eq:intro:central_cone} arises as a morphism of monads and make precise its relationship with Hochschild cohomology of $\catDMod{\stack X}$.
As a simple application of this formalism, we show that $\HCoh\bigl(\catDMod{\stack X}\bigr)$ always has a subalgebra isomorphic to $\opalg{\GammadR(\stack X,\ k_{\stack X})}$.

In Section~\ref{sec:base-change} we give a criterion for morphism~\eqref{eq:intro:central_cone} to be an isomorphism.
In particular, we will use the contraction principle of Drinfeld--Gaitsgory \cite[Section~5.1]{DrinfeldGaitsgory:2015:CompactGenerationOfDModOnBunG} to understand the case of so-called contractive stacks.
Using the theory developed in Section~\ref{sec:base-change}, the proof of Theorem~\ref{thm:main} is given in Section~\ref{sec:torus}.

Finally, in Section~\ref{sec:compactification} we describe the cone of morphism~\eqref{eq:intro:central_cone} as sections of a sheaf on a \enquote{compactified} loop space of $\stack X$ and consider some of its properties.

\subsection*{Future directions}

The present paper represents only the first step in the study of Hochschild cohomology and associated support theories for D-modules.
In future work we intend to expand upon the foundations laid out in this paper.
Of particular interest are quotients stacks $X/G$ with $G$ an affine algebraic group and the stack $\Bun_G(C)$ of $G$-bundles on a curve $C$.
Following the constructions laid out in Section~\ref{sec:compactification:quotient}, for the former we need a good $G\times G$-equivariant compactification of $G$.
Thus groups with a wonderful compactification are of particular interest for study with our methods.
For the latter, relative compactifications of the diagonal of $\Bun_G(C)$ have been proposed by Lafforgue and Drinfeld.
These spaces are under active study and we expect that any computations of Hochschild cohomology to have deep connections with the \enquote{pseudo-identity} that is needed in the formulation of the geometric (categorical) Langlands conjecture.

\subsection*{Acknowledgments}
I would like to thank David Nadler for the original motivation for this project and many discussions concerning it.
I would also like to thank David Ben-Zvi, Dennis Gaitsgory and Nathan Ilten for discussions concerning this project.
This work has been partially supported by the National Science Foundation under Grant No.~1638352, as well as the Giorgio and Elena Petronio Fellowship Fund.

\section{D-modules on stacks}\label{sec:d-mods}

We fix an algebraically closed base field $k$ of characteristic $0$.
All stacks will be assumed to be algebraic over $k$ and quasi-compact with affine automorphism groups (QCA).
Thus by definition for any stack $\stack X$ we have:
\begin{itemize}
    \item The diagonal morphism $\Delta \colon \stack X \to  \stack X \times  \stack X$ is schematic.
    \item There exists a scheme $Z$ with a smooth and surjective map $Z \to  \stack X$.
    \item $\stack X$ is quasi-compact.
    \item The automorphism groups of the geometric points of $\stack X$ are affine.
    \item The loop space (or inertia stack) $\ls \stack X = \stack X \times _{\stack X \times  \stack X} \stack X$ is of finite presentation over $\stack X$.
\end{itemize}
The QCA condition guarantees a reasonable theory of D-modules on $\stack X$.
For details on QCA stacks we refer to~\cite{DrinfeldGaitsgory:2013:FinitenessQuestions}.
Every quotient of a scheme of finite type over $k$ by an affine algebraic group is a QCA stack, and we will be mainly interested in these.

In order to correctly define categories of D-modules on stacks it is necessary to work with dg-categories.
We refer to \cite{Keller:2006:OnDGCategories} for an introduction to dg categories.
It is often convenient to regard pretriangulated dg categories as $k$-linear stable $(\infty ,1)$-categories \cite{Lurie:2009:HigherToposTheory,Lurie:2017-draft:HigherAlgebra}, which can be done via the nerve construction. 
We will switch between those two languages without explicitly mentioning the intervening constructions and apply results from \cite{Lurie:2017-draft:HigherAlgebra} to dg categories.
Fortunately, a superficial knowledge of dg/$\infty $-categories should be sufficient for reading this article.

The category of D-modules on a stack $\stack X$ can be either constructed via descent \cite{BeilinsonDrifeld:unpublished:Hitchin,DrinfeldGaitsgory:2013:FinitenessQuestions} or equivalently as ind-coherent sheaves on the de Rham space of $\stack X$ \cite{GaitsgoryRozenblyum:2014:CrystalsAndDModules}.
While the first construction is more \enquote{hands on}, the second construction is often more useful from a theoretical point of view.
It is explained in detail in the book \cite{GaitsgoryRozenblyum:2017:StudyInDAG:2} (see also \cite{FrancisGaitsgory:2012:ChiralKoszulDuality} for an overview).
Many basic properties of the category $\catDMod{\stack X}$ are explored in \cite{DrinfeldGaitsgory:2013:FinitenessQuestions} and unless stated otherwise proofs for the assertions in this section can be found there.

For any morphism $f\colon \stack X \to  \stack Y$ we have a continuous functor $f^!\colon \catDMod{\stack Y} \to  \catDMod{\stack X}$ and a (not necessarily continuous) functor $f_*\colon \catDMod{\stack X} \to  \catDMod{\stack Y}$.
If $p\colon \stack X \to  \pt$ is the structure map, we set
\[
    \GammadR(\stack X,\, {-}) = p_*({-}) \colon \catDMod{\stack X} \to  \catVect.
\]
The functor $\GammadR$ is representable by a D-module $k_{\stack X}$, i.e.
\[
    \GammadR(\stack X,\, {-}) = \Hom_{\catDMod{\stack X}}(k_{\stack X},\, {-}).
\]
As $\GammadR(\stack X,\, {-})$ is usually not continuous the object $k_{\stack X}$ is usually not a compact object of $\catDMod{\stack X}$.

Let $\Delta \colon \stack X \to  \stack X \times  \stack X$ be the diagonal.
The category $\catDMod{\stack X}$ has a monoidal structure given by the tensor product $\sheaf F \otimes \sheaf G \coloneqq \Delta ^!\bigl( \sheaf F \boxtimes \sheaf G \bigr)$ with unit $\omega _{\stack X} = p^! k$.

We will be mainly concerned with the subcategory of holonomic D-modules since they enjoy extended functoriality.
\begin{definition}
    A $D$-module $\sheaf F \in  \catDMod{\stack X}$ is called \emph{holonomic} if $f^!\sheaf F$ is holonomic for any smooth morphism $f\colon Z \to  \stack X$ from a scheme $Z$.
    The full subcategory of holonomic D-modules will be denoted $\catDModHol{\stack X}$.
\end{definition}

The following assertions mostly follow from their corresponding counterparts for schemes.
We refer to \cite{Braverman:LecturesOnAlgebraicDmodules} for proofs in the case of non-smooth schemes.

\begin{proposition}
    Let $f\colon \stack X \to  \stack Y$ be a schematic morphism.
    Then $f^!$ and $f_*$ restrict to functors
    \[
        f^!\colon \catDModHol{\stack Y} \to  \catDModHol{\stack X}
        \quad\text{and}\quad
        f_*\colon \catDModHol{\stack X} \to  \catDModHol{\stack Y}.
    \]
\end{proposition}

The Verdier duality functor on schemes induces an involutive anti-auto-equivalence
\[
    \dualize_{\stack X}\colon \catDModHol{\stack X}^\mathrm{op} \to  \catDModHol{\stack X}
\]
such that for each smooth morphism $Z \to  \stack X$ of relative dimension $d$ from a scheme $Z$ one has
\[
    f^! \circ  \dualize_{\stack X} \cong \dualize_{Z} \circ  f^![-2d].
\]
The Verdier duality functor then allows us to define the \emph{non-standard functors} $f_!$ and $f^*$ for any schematic morphism $f\colon \stack X \to  \stack Y$ by
\begin{align*}
    f^* & = \dualize_{\stack X} \circ  f^! \circ  \dualize_{\stack Y} \colon \catDModHol{\stack Y} \to  \catDModHol{\stack X} \\
    \intertext{and}
    f_! & = \dualize_{\stack Y} \circ  f_* \circ  \dualize_{\stack X} \colon \catDModHol{\stack X} \to  \catDModHol{\stack Y}.
\end{align*}
We obtain adjoint pairs $(f_!,\, f^!)$ and $(f^*,\, f_*)$.
In some situations we can identify the non-standard functors with their standard counterparts.
If $f$ is smooth of relative dimension $d$ then $f^* = f^![-2d]$.
If $f$ is proper then $f_! = f_*$ and in particular $f_*$ is left adjoint to $f^!$.
The objects $\omega _{\stack X}$ and $k_{\stack X}$ are always holonomic and
\[
    \dualize_{\stack X} \omega _{\stack X} = k_{\stack X}.
\]
We have $k_{\stack X} = f^* k_{\stack Y}$ and if $\stack X$ is smooth, then $k_{\stack X} = \omega _{\stack X}[-2\dim \stack X]$.

We will make use of the following lemma which follows from \cite[Lemma~5.1.6]{DrinfeldGaitsgory:2013:FinitenessQuestions}.

\begin{lemma}
    For a smooth and schematic morphism $f$ the functor $f^!$ is conservative.
\end{lemma}

\begin{proposition}[{\cite[Section~4.2.1.3]{GaitsgoryRozenblyum:2017:StudyInDAG:2}}]
    \label{prop:pre:base-change}%
    Consider a Cartesian square
    \[
        \begin{tikzcd}
            \stack Z \arrow{d}{p} \arrow{r}{q} & \stack X_1 \arrow{d}{f} \\
            \stack X_2 \arrow{r}{g} & \stack Y
        \end{tikzcd}
    \]
    with schematic morphism $f$ (and hence also schematic $p$).
    Then there is a base change isomorphism
    \[
        p_* q^! \isoto g^! f_*
    \]
    of functors from $\catDMod{\stack{X_1}}$ to $\catDMod{\stack{X_2}}$.
    If furthermore $f$ (and hence $p$) is proper, then this isomorphism coincides with the natural transformation
    \[
        p_* q^! \to 
        p_* q^! f^! f_* =
        p_* p^! g^! f_* \to 
        g^! f_*
    \]
    induced by $(f_*,\,f^!)$ and $(p_*,\, p^!)$ adjunctions.
    Similarly, if $f$ is an open immersion, then the base change isomorphism coincides with the natural transformation resulting from the adjunctions $(f^!,f_*)$ and $(p^!, p_*)$.
\end{proposition}

\begin{proposition}
    \label{prop:pre:projection-formula}%
    If $f\colon \stack X \to  \stack Y$ is a schematic morphism then the projection formula holds, i.e.~there is a functorial isomorphism
    \[
        \sheaf F \otimes f_*(\sheaf G) \cong f_*\bigl( f^! \sheaf F \otimes \sheaf G)
    \]
    for $\sheaf F \in  \catDMod{\stack Y}$ and $\sheaf G \in  \catDMod{\stack X}$.
\end{proposition}

\begin{remark}
    Propositions \ref{prop:pre:base-change} and \ref{prop:pre:projection-formula} hold more generally when $f$ is merely a \enquote{safe} morphism.
    Alternatively they hold in full generality after replacing $f_*$ by the \enquote{renormalized de Rham pushforward}.
    We will not use either notion here and refer the interested reader to \cite{DrinfeldGaitsgory:2013:FinitenessQuestions}.
\end{remark}

For D-modules on stacks we have the usual recollement package.
Let $i\colon \stack Z \hookrightarrow \stack X$ be a closed embedding and $j\colon \stack U \hookrightarrow \stack X$ the complementary open.
We have adjoint pairs $(i_*,\, i^!)$ and $(j^!,\, j_*)$.

\begin{proposition}[{\cite[Section~2.5]{GaitsgoryRozenblyum:2014:CrystalsAndDModules}}]
    \label{prop:recollement-std}%
    There is an exact triangle of functors
    \[
        i_*i^! \to  \id \to  j_*j^!
    \]
    on $\catDMod{\stack X}$, the adjunction morphisms
    \[
        \id \to  i^!i_*
        \quad\text{and}\quad
        j^!j_* \to  \id
    \]
    are isomorphisms, the functors $j^!i_*$ and $i^!j_*$ vanish and $i_*$ and $j_*$ are full embeddings.
\end{proposition}

On holonomic D-modules we have the additional adjoint pairs $(i^*,\, i_*)$ and $(j_!,\, j^*)$.
By applying duality to Proposition~\ref{prop:recollement-std} we obtain the distinguished triangle
\[
    j_! j^* \to  \id \to  i_*i^*
\]
and the identity $i^*j_! = 0$ on holonomic D-modules.
Further, the functor $j_!$ is a full embedding $\catDModHol{\stack U} \hookrightarrow \catDModHol{\stack X}$.

\section{Monads}\label{sec:monads}

We will deduce Theorem~\ref{thm:main} from a isomorphism of monads on $\catDMod{\stack X}$.
In this section we give a short introduction to the theory of monads and the specific constructions that we will use.
However, in the interest of readability we will mainly do so informally, skipping over the intricacies of $\infty $-categories.
The interested reader can find the correct $\infty $-categorical formulations in the given references.

Thus we think of a \emph{monad} on a category $\cat C$ as consisting of a triple $(T, \eta , \mu )$, where $T\colon \cat C \to  \cat C$ is an endofunctor of $\cat C$, and $\eta \colon \id_{\cat C} \to  T$ and $\mu \colon T\circ T \to  T$ are natural transformations such that the diagrams
\begin{equation}
    \label{eq:pre:monad-identities}
    \begin{tikzcd}
        T^3 \arrow{r}{T\mu } \arrow{d}[swap]{\mu T} & T^2 \arrow{d}{\mu } \\
        T^2 \arrow{r}[swap]{\mu } & T
    \end{tikzcd}
    \quad\text{and}\quad
    \begin{tikzcd}
        T \arrow{r}{T\eta } \arrow{d}[swap]{\eta T} \arrow{dr}{\mathrm{id}} & T^2 \arrow{d}{\mu } \\
        T^2 \arrow{r}[swap]{\mu } & T
    \end{tikzcd}
\end{equation}
commute.
Alternatively, we can think of $T$ being a monoid in the category of endofunctors of $\cat C$ with the monoidal structure given by composition of endofunctors.
This definition also gives the correct generalization to $\infty $-categories \cite[Definition~4.7.0.1]{Lurie:2017-draft:HigherAlgebra}.

Let $X$ be an object of the category $\cat C$.
Then $T$ gives the dg vector space $\Hom_{\cat C}(X, TX)$ the structure of a dg algebra with multiplication map
\[
    (f,g) \mapsto \mu _X \circ  Tf \circ  g,
\]
\[
    \begin{tikzcd}
        X \arrow{r}{g} & TX \arrow{r}{Tf} & T^2X \arrow{r}{\mu _X} & TX.
    \end{tikzcd}
\]
The identities~\eqref{eq:pre:monad-identities} ensure that the algebra is associative and unital.

The most common source of monads is a pair of adjoint functors $F\colon \cat C \rightleftarrows \cat D \cocolon G$.
One simply sets $T = G \circ  F$ and $\eta $ and $\mu $ are given by the adjunction morphisms
\[
    \id[\cat C] \to  G \circ  F = T
    \qquad\text{and}\qquad
    T^2 = G \circ  (F \circ  G) \circ  F \to  G \circ  F = T.
\]
We note that the correct construction is more complicated in the $\infty $-categorical case and refer to \cite[Section~4.7]{Lurie:2017-draft:HigherAlgebra}.
For any $X \in  \cat C$ the algebra construction above gives an isomorphism of algebras
\[
    \Hom_{\cat C}(X, (GF)(X)) \cong
    \Hom_{\cat D}(FX, FX).
\]

Another common way to obtain monads in geometry is via a groupoid.
Recall that a groupoid $G_{\cx}$ in stacks consists of a stack $\stack G_0$ of \enquote{objects} and a stack $\stack G_1$ of \enquote{morphisms} together with
\begin{itemize}
    \item \emph{source} and \emph{target} maps $s,t\colon \stack G_1 \rightrightarrows \stack G_0$,
    \item a \emph{unit} $e\colon \stack G_0 \to  \stack G_1$,
    \item a \emph{multiplication} (or \emph{composition}) map $m\colon \stack G_1 \mathop{\times }\limits_{s,\stack G_0,t} \stack G_1 \to  \stack G_1$,
    \item an \emph{inverse} map $\iota \colon \stack G_1 \to  \stack G_1$,
\end{itemize}
such that
\begin{itemize}
    \item $s \circ  e = t \circ  e = \id_{\stack G_0}$,
    \item $s \circ  m = s \circ  p_2$ and $t \circ  m = t \circ  p_1$ (where $p_i\colon \stack G_1 \times _{s,\stack G_0,t} \stack G_1 \to  \stack G_1$ are the projection maps).
    \item $m$ is associative,
    \item $\iota $ interchanges $s$ and $t$ and is an inverse for $m$,
\end{itemize}
where all identities have to be understood in the correct $\infty $-categorical way \cite[Section~6.1.2]{Lurie:2009:HigherToposTheory}.

\begin{example}
    For our purpose the most important example is the following:
    Let $f\colon \stack X \to  \stack S$ be a morphism of stacks.
    We set $\stack G_0 = \stack X$ and $\stack G_1 = \stack X \times _\stack S \stack X$.
    The source and target maps are given by $p_1$ and $p_2$, the unit by the diagonal $\Delta \colon \stack X \to  \stack X\times _\stack S\stack X$, the inverse by interchanging the factors and multiplication is $p_{13}\colon \stack X \times _\stack S \stack X \times _\stack S \stack X \to  \stack X\times _\stack S\stack X$.
\end{example}

Let us for the moment assume that $s$ (and hence $p_2$) is proper and schematic.
In this case the maps $e$ and $m$ are also proper, since $s \circ  e = \id_{\stack G_0}$ and $s \circ  m = s \circ  p_2$ are proper.
In particular the functors $e^!\colon \catDMod{\stack G_1} \to  \catDMod{\stack G_0}$ and $m^!\colon \catDMod{\stack G_1} \to  \catDMod{\stack G_1 \times _{\stack G_0} \stack G_1}$ have left adjoints given by $e_*$ and $m_*$ respectively.
This allows us to give the endofunctor $T = s_*t^!$ of $\catDMod{\stack G_0}$ the structure of a monad in the following way:

\begin{itemize}
    \item By $(e_*,e^!)$-adjunction we have a transformation
        \[
            \id = (s\circ e)_*(t\circ e)^! = s_*e_*e^!t^! \to  s_*t^! = T.
        \]
    \item Consider the following commutative diagram
        \[
            \begin{tikzcd}[column sep=small]
                {} & & \stack G_1 \arrow[bend right]{dddll}[swap]{s} \arrow[bend left]{dddrr}{t} & & \\
                   & & \stack G_1 \times _{\stack G_0} \stack G_1 \arrow{u}{m} \arrow{dl}{p_2} \arrow{dr}{p_1} & & \\
                   & \stack G_1 \arrow{dl}{s} \arrow{dr}{t} & & \stack G_1 \arrow{dl}{s} \arrow{dr}{t} & \\
                \stack G_0 & & \stack G_0 & & \stack G_0
            \end{tikzcd}
        \]
        with Cartesian middle square.
        Proper base change and $(m_*,m^!)$-adjunction gives a transformation
        \[
            T^2 =
            s_*t^!s_*t^! =
            (s\circ p_2)_*(t\circ p_1)^! =
            (s\circ m)_*(t\circ m)^! =
            s_*m_*m^!t^! \to 
            s_*t^! =
            T.
        \]
\end{itemize}

In the non-$\infty $-categorical setting one could easily check by hand that this is indeed a monad.
To obtain the corresponding derived statement one applies an argument similar to \cite[Section~4.7.2]{GaitsgoryRozenblyum:2017:StudyInDAG:1}.

Let now $f\colon \stack X \to  \stack Y$ be schematic and proper.
The Cartesian diagram
\[
    \begin{tikzcd}
        \stack X \times _{\stack Y} \stack X \arrow{r}{p_s} \arrow{d}{p_t} & \stack X \arrow{d}{f} \\
        \stack X \arrow{r}{f} & \stack Y
    \end{tikzcd}
\]
induces a groupoid with $\stack{G}_0 = \stack X$ and $\stack{G}_1 = \stack X \times _{\stack Y} \stack X$.
The above constructions give two monads on $\catDMod{\stack{X}}$: one by $(f_*,f^!)$ adjunction and one from the groupoid structure.
The base change isomorphism
\[
    p_{t,*} p_s^! \to  f^! f_*
\]
gives an identification of these monads and hence of the algebras that they induce, i.e.~for any $\sheaf F \in  \catDMod{\stack X}$ we have
\[
    \Hom(\sheaf F,\, p_{t,*} p_s^!\sheaf F) \cong
    \Hom(\sheaf F,\, f^! f_* \sheaf F) \cong
    \Hom(f_*\sheaf F,\, f_* \sheaf F).
\]

We will need to apply this construction for non-proper $f$.
Unfortunately, in this case none of the adjunctions used to define the monads are available.
We rectify this by restricting to the full subcategory of of holonomic D-modules and using the $!$-pushforward functors instead of the $*$-pushforward ones.
Of course, by doing so we do not automatically have base change isomorphisms available anymore.
Thus we have to explicitly require the \emph{Beck--Chevalley condition}, saying that all necessary base changes hold.

The Beck--Chevalley condition is formulated using the \emph{nerve} of the groupoid $\stack G_\cx$.
Intuitively, this is the simplicial stack, usually also denoted $\stack{G}_\cx$, with
\[
    \stack G_i = \underbrace{\stack G_1 \times _{\stack G_0} \dotsb \times _{\stack G_0} \stack G_1}_{\text{$i$ factors}}.
\]
We refer to \cite[Section~6.1.2]{Lurie:2009:HigherToposTheory} for the correct $\infty $-categorical setup.
The following lemma is now an immediate corollary of \cite[Lemma~4.7.1.4]{GaitsgoryRozenblyum:2017:StudyInDAG:2} or \cite[Theorem~4.7.5.2]{Lurie:2017-draft:HigherAlgebra}.

\begin{lemma}
    \label{lem:pre:groupoid_monad_hol}%
    Let $f\colon \stack X \to  \stack Y$ be a schematic morphism of stacks and let $\stack G_\cx$ be the corresponding groupoid.
    For each map $F\colon [n] \to  [m]$ in $\cat{\Delta }^{\mathrm{op}}$ consider the corresponding square
    \[
        \begin{tikzcd}
            \stack G_{n+1} \arrow{r}{p_s} \arrow{d}{p_{F+1}} & \stack G_n \arrow{d}{p_F} \\
            \stack G_{m+1} \arrow{r}{p_s} & \stack G_m
        \end{tikzcd}
    \]
    where the vertical arrows are induced by $F$.
    Assume that for each such square the base change composition
    \[
        p_{F+1,!} p_s^! \to 
        p_{F+1,!} p_s^! p_F^! p_{F,!} =
        p_{F+1,!} p_{F+1}^! p_s^!  p_{F,!} \to 
        p_s^! p_{F,!}
    \]
    given by the adjunction morphisms is an isomorphism of functors from $\catDModHol{\stack G_n}$ to $\catDModHol{\stack G_{m+1}}$.
    Further assume that the same is true along 
    \[
        \begin{tikzcd}
            \stack G_1 \arrow{r}{p_s} \arrow{d}{p_t} & \stack G_0 \arrow{d}{f} \\
            \stack G_0 \arrow{r}{f} & \stack Y
        \end{tikzcd}
    \]
    Then the endofunctor $p_{t,!} p_s^!$ of $\catDModHol{\stack X}$ has a canonical structure of a monad and as such is isomorphic to the adjunction monad $f^!f_!$.
\end{lemma}

\section{Hochschild cohomology and kernels}\label{sec:HH}

We recall that the Hochschild cohomology of a dg category $\cat C$ is the algebra of derived endomorphisms of the identity functor,
\[
    \HCoh(\cat C) = \mathbb{R}\mathrm{Hom}(\id[\cat C],\, \id[\cat C]).
\]
For the exact definition of the category $\mathbb{R}\mathrm{Hom} = \cat{Funct}(\cat C, \cat C)$ we refer to \cite{Keller:2006:OnDGCategories}.
Instead we will give a more concrete construction via kernels which can be applied to $\catDMod{\stack X}$.
For this let us restrict to the case of co-complete dg categories and let $\cat{Funct}_{\mathrm{cont}}(\cat C, \cat C)$ be the full subcategory of $\cat{Funct}(\cat C, \cat C)$ spanned by the continuous functors.
Then, since $\id[\cat C]$ is evidently continuous, we have
\[
    \HCoh(\cat C) =
    \Hom_{\cat{Funct}_{\mathrm{cont}}(\cat C, \cat C)}(\id[\cat C], \id[\cat C]).
\]
Let us further assume that $\cat C$ is dualizable with dual $\cat C^\dual$.
Thus there is a unit map
\[
    \eta \colon \catVect \to  \cat C^\dual \otimes \cat C
\]
and a counit map
\[
    \epsilon \colon \cat C \otimes \cat C^\dual \to  \catVect
\]
satisfying the usual compatibilities (cf.~\cite[Section~2]{BenZviNadler:arXiv:NonlinearTraces}).
Let $u = \eta (k)$.
To each continuous endofunctor $F$ of $\cat C$ we can associate its kernel $\bigl(\id[\cat C^\dual] \otimes F\bigr)(u) \in  \cat C^\dual \otimes \cat C$ and conversely to each kernel $Q \in  \cat C^\dual \otimes \cat C$ we can associate the endofunctor
\[
    \cat C
    \xrightarrow{\id[\cat C] \otimes Q}
    \cat C \otimes \cat C^\dual \otimes \cat C
    \xrightarrow{\epsilon  \otimes \id[\cat C]}
    \cat C.
\]
These assignments are mutually inverse and give an equivalence of dg categories
\[
    \cat{Funct}_{\mathrm{cont}}(\cat C, \cat C)
    \cong
    \cat C^\dual \otimes \cat C.
\]
In particular, the kernel for the identity is $u$ and hence we have
\[
    \HCoh(\cat C) =
    \Hom_{\cat C^\dual \otimes \cat C}(u, u).
\]
Let us now consider the case of $\cat C = \catDMod{\stack X}$ for a stack $\stack X$.
Let $p\colon \stack X \to  \pt$ be the structure morphism and $\Delta \colon \stack X \to  \stack X \times  \stack X$ the diagonal.
By \cite[Section~8.4]{DrinfeldGaitsgory:2013:FinitenessQuestions} the category $\catDMod{\stack X}$ is dualizable and there is a canonical identification
\[
    \catDMod{\stack X}^\dual \otimes \catDMod{\stack X} \cong \catDMod{\stack X \times  \stack X}
\]
such that the unit map is given by $\Delta _*p^!$, i.e.~we have
\[
    u = \Delta _*\omega _{\stack X}.
\]
We summarize the above discussion in the following lemma.

\begin{lemma}\label{lem:pre:hcoh}
    Let $\stack X$ be a stack.
    Then the Hochschild cohomology of $\catDMod{\stack X}$ is given by the dg algebra
    \[
        \HCoh(\catDMod{\stack X}) =
        \Hom_{\catDMod{\stack X \times  \stack X}}(\Delta _*\omega _{\stack X},\, \Delta _*\omega _{\stack X}).
    \]
\end{lemma}

Recall now that we assume $\Delta $ to be schematic and that the dualizing module $\omega _{\stack X}$ is always holonomic.
Thus we have $\dualize \Delta _* \omega _{\stack X} = \Delta _! k_{\stack X}$.
With this we observe that
\begin{align*}
    \HCoh\bigl(\catDMod{\stack X}\bigr)
    & = \Hom_{\catDMod{\stack X \times  \stack X}}(\Delta _* \omega _{\stack X},\, \Delta _* \omega _{\stack X}) & &\text{(Lemma~\ref{lem:pre:hcoh})} \\
    & = \opalg{\Hom_{\catDMod{\stack X \times  \stack X}}(\Delta _! k_{\stack X},\, \Delta _! k_{\stack X})} & & \text{(duality)} \\
    & = \opalg{\Hom_{\catDMod{\stack X \times  \stack X}}(k_{\stack X},\, \Delta ^! \Delta _! k_{\stack X})} & & \text{(adjunction)} \\
    & = \opalg{\GammadR\bigl(\stack X,\, \Delta ^! \Delta _! k_{\stack X}\bigr)},
\end{align*}
where the algebra structure on $\GammadR\bigl(\stack X,\, \Delta ^! \Delta _! k_{\stack X}\bigr) = \Hom_{\catDMod{\stack X \times  \stack X}}(k_{\stack X},\, \Delta ^! \Delta _! k_{\stack X})$ is the one induced by the $(\Delta _!,\Delta ^!)$-adjunction monad.
Consider the Cartesian square
\[
    \begin{tikzcd}
        \ls{\stack X} \arrow{r}{p_2} \arrow{d}{p_1} & \stack X \arrow{d}{\Delta } \\
        \stack X \arrow{r}{\Delta } & \stack X \times  \stack X
    \end{tikzcd}
\]
Let us assume for the moment that $\Delta $ (and hence $p_i$) is proper.
Then $\Delta _* = \Delta _!$ and $p_{1,*} = p_{1,!}$ and by Section~\ref{sec:monads} we have an isomorphism of monads on $\catDModHol{\stack X}$
\begin{equation}
    \label{eq:central_iso}
    p_{1,!} p_2^! \to  \Delta ^!\Delta _!,
\end{equation}
which induces an isomorphism of algebras
\[
    \GammadR\bigl(\stack X,\, p_{1,!} p_2^! k_{\stack X}\bigr)
    \to 
    \GammadR\bigl(\stack X,\, \Delta ^! \Delta _! k_{\stack X}\bigr).
\]
Of course, if $\stack X$ is QCA but not an algebraic space, then $\Delta $ is not proper (nor is it in general smooth).
Thus in general \eqref{eq:central_iso} is not an isomorphism and there is no canonical structure of monad on $p_{1,!} p_2^!$.
We would like to apply Lemma~\ref{lem:pre:groupoid_monad_hol} to construct a monad in special cases.
Thus the goal of the next section is to give a criterion for the assumptions of Lemma~\ref{lem:pre:groupoid_monad_hol}, i.e.~for base change to hold.

\begin{example}
    The base change morphism \eqref{eq:central_iso} is also an isomorphism if $\Delta $ is smooth.
    In particular this implies that the \enquote{naive expectation} holds for $\stack X = \cs{G}$ for any algebraic group $G$, i.e.~we have
    \[
        \HCoh(\catDMod{\cs G}) = \opalg{\GammadR\bigl(\cs G,\, p_{1,!} p_2^! k_{\cs G}\bigr)}.
    \]
    An argument similar to \cite{BenZvi:mathoverflow:CohomologyOfGmodG} shows that there is a further isomorphism
    \begin{equation}\label{eq:HH_BG_decomp}
        \GammadR\bigl(\cs G,\, p_{1,!} p_2^! k_{\stack X}\bigr) \cong
        \GammadR(G, k_G)^\dual \otimes \GammadR(\cs G, k_{\cs G}).
    \end{equation}
    Alternatively, we can use the identification
    \[
        \catDMod{\cs G} \cong \catModules{\GammadR(G, k_{G})^\dual},
    \]
    where the algebra structure on $\GammadR(G, k_G)^\dual$ is induced by the group multiplication \cite[Section~7.2]{DrinfeldGaitsgory:2013:FinitenessQuestions}.
    If $G$ is reductive, then $\GammadR(G, k_G)^\dual$ is an exterior algebra and thus its Hochschild cohomology can be computed directly.
\end{example}

From \eqref{eq:HH_BG_decomp} we see that $\HCoh(\catDMod{\cs G})$ contains a copy of $\opalg{\GammadR(\cs G, k_{\cs G})}$.
Indeed this is a general phenomenon.

\begin{proposition}
    \label{prop:mono-from-cohomology}%
    The Hochschild cohomology $\HCoh\bigl(\catDMod{\stack X}\bigr)$ has $\opalg{\GammadR\bigl(\stack X,\, k_{\stack X}\bigr)}$ as a direct summand.
\end{proposition}

\begin{proof}
    Consider the sequence of maps
    \[
        \stack X \xrightarrow{\Delta } \stack X \times  \stack X \xrightarrow{\proj{1}} \stack X.
    \]
    Since the composition is the identity, the functor $(\proj{1})_! \circ  \Delta _!$ is faithful.
    Hence $\Delta _!$ is faithful too.
    Thus the adjunction map
    \[
        \eta\colon \id[\catDModHol{\stack X}] \to  \Delta ^!\Delta _!
    \]
    is a monomorphism.
    Completing to a distinguished triangle
    \[
        k_{\stack X} \xrightarrow{\eta} \Delta ^!\Delta _! k_{\stack X} \to  \sheaf F \xrightarrow{\delta }
    \]
    for some $\sheaf F \in \catDModHol{\stack X}$, we see that the connection morphism $\delta $ is the zero morphism (indeed $\eta \circ \delta[1] = 0$ and $\eta$ is mono).
    Thus the same is true after applying $\Hom(k_{\stack X}, \{-\})$.
    In particular, the induced map
    \[
        \GammadR\bigl(\stack X,\, k_{\stack X}\bigr)
        \to 
        \GammadR\bigl(\stack X,\, \Delta ^!\Delta _! k_{\stack X}\bigr)
        \cong
        \opalg{\HCoh\bigl(\catDMod{\stack X}\bigr)}
    \]
    is a split monomorphism.
\end{proof}

\section{Base change along the diagonal}\label{sec:base-change}

Consider a Cartesian diagram of stacks
\begin{equation*}
    \begin{tikzcd}
        \stack Z \arrow{d}{p} \arrow{r}{q} & \stack X_1 \arrow{d}{f} \\
        \stack X_2 \arrow{r}{g} & \stack Y
    \end{tikzcd}
\end{equation*}
with $f$ and $g$ schematic.
We have a morphism of functors $\catDModHol{\stack X_1} \to  \catDModHol{\stack X_2}$,
\begin{equation}
    \label{eq:base-change-morphism}
     p_! q^! \to  g^! f_!
\end{equation}
induced by adjunctions
\begin{equation*}
    p_! q^! \to 
    p_! q^! f^! f_! =
    p_! p^! g^! f_! \to 
    g^! f_!.
\end{equation*}
If $f$ is proper, then \eqref{eq:base-change-morphism} is an isomorphism by Proposition~\ref{prop:pre:base-change}.

\subsection{The required base changes}

In order for the isomorphism \eqref{eq:intro:naive} to hold, we only need know that the base change morphism is an isomorphism on $k_{\stack X}$.
However, it is often useful to pass to substacks, so that we need to prove the isomorphism of a slightly larger class of D-modules.
Thus we introduce the following subcategories.

\begin{definition}
  Let $\catK{\stack X}$ be the full subcategory of $\catDMod{\stack X}$ generated by $k_{\stack X}$ and let $\catKHol{\stack X}$ be its full subcategory consisting of holonomic objects.
\end{definition}

With this we can write down the two properties that we need to check in order for the naive expectation to hold.

\begin{definition}
    Let $\stack X$ be any QCA stack.
    Consider the Cartesian diagram
    \[
        \begin{tikzcd}
            \ls \stack X \arrow{d}{p_1} \arrow{r}{p_2} & \stack X \arrow{d}{\Delta} \\
            \stack X \arrow{r}{\Delta} & \stack X \times \stack X
        \end{tikzcd}
    \]
    We say that $\stack X$ has \emph{property \bc} if the canonical morphism
    \[
        p_{1,!}p_2^! \to \Delta^!\Delta_!
    \]
    of functors $\catKHol{\stack X} \to \catDModHol{\stack X}$ is an equivalence.
\end{definition}

Let us write $\stack G_\bullet = \stack G_\bullet(\stack X)$ for the nerve of the groupoid associated to $\ls\stack X \rightrightarrows \stack X$.
Thus
\[
    \stack G_m(\stack X) = \underbrace{\stack X \XXtimes \dotsb \XXtimes \stack X}_{m+1\text{ times}}.
\]

\begin{definition}
    Let $\stack X$ be a QCA stack and consider the Cartesian square
    \[
        \begin{tikzcd}
            \stack G_{n+1} \arrow{r}{p_{s}} \arrow{d}{p_{F+1}} & \stack G_n \arrow{d}{p_F} \\
            \stack G_{m+1} \arrow{r}{p_s} & \stack G_m
        \end{tikzcd}
    \]
    where the horizontal maps are induced by a map $F\colon [n] \to [m]$ in $\catCoSimplicial$.
    We say that $\stack X$ has \emph{property \hbc} if for each such $F$ the canonical morphism
    \[
        p_{F+1,!}p_s^! \to p_s^!p_{F,!}
    \]
    of functors $\catKHol{\stack G_n} \to \catDModHol{\stack G_{m+1}}$ is an equivalence.
\end{definition}

\begin{corollary}\label{cor:properties_imply_monads}
  If $\stack X$ has properties \bc\ and \hbc, then there exists a canonical structure of monad on $p_{1,!}p_2^!$ and the morphism $p_{2,!}p_1^! \to  \Delta ^!\Delta _!$ is an isomorphism of monads.
  In particular there is an isomorphism of algebras
  \[
      \HCoh\bigl(\catDMod{\stack X}\bigr)
      \cong
      \opalg{\GammadR\bigl(\stack X,\, p_{1,!} p_2^! k_{\stack X}\bigr)}.
  \]
\end{corollary}

\begin{proof}
  Since $\Hom(k_{\stack X},\, \Delta^!\Delta_!k_{\stack x}) = \Hom(\Delta_!k_{\stack X},\, \Delta_!k_{\stack X}) \ne 0$, we see that $\Delta^!\Delta_!$ restricts to a monad on $\catKHol{\stack X}$.
  Thus by the assumption $p_{1,!} p_2^!$ is also an endofunctor of $\catKHol{\stack X}$.
  The statement now follows from a restriction of Lemma~\ref{lem:pre:groupoid_monad_hol} applied to the construction in Section~\ref{sec:HH}.
\end{proof}

\begin{remark}
    We expect that most quotient stacks do not have property \bc.
    For example, a direct computation shows that the transformation of property \bc\ is not an isomorphism for the D-module $k$ on the stack $\ps1/\Ga$.
\end{remark}

\begin{lemma}\label{lem:smooth_or_proper}
    Let $\stack X$ be a QCA stack with diagonal morphism $\Delta\colon \stack X \to \stack X \times \stack X$.
    \begin{enumerate}
        \item\label{lem:smooth_or_proper:bc} If $\Delta$ is either proper of smooth, then $\stack X$ has property \bc.
        \item\label{lem:smooth_or_proper:hbc} If $\Delta$ is proper, then $\stack X$ has property \hbc.
        \item\label{lem:smooth_or_proper:scheme} If $\stack X$ is a separated scheme, then $\stack X$ has both properties \bc\ and \hbc.
        \item\label{lem:smooth_or_proper:bg} Let $G$ be an affine algebraic group. Then the classifying stack $\cs{G}$ has both properties \bc\ and \hbc.
    \end{enumerate}
\end{lemma}

\begin{proof}
    Statement \ref{lem:smooth_or_proper:bc} is a direct consequence of the base change theorem, Proposition~\ref{prop:pre:base-change}.
    To be more explicit, the smooth version follows from Proposition~\ref{prop:pre:base-change} by duality and the fact that $f^!$ agrees with $f^*$ up to shift.

    We only need to show \ref{lem:smooth_or_proper:hbc} under the additional assumption that $F$ is a face or degeneracy map.
    In either case it follows immediately that $p_{F}$ is proper and we can again apply base change.

    Statement \ref{lem:smooth_or_proper:scheme} is just a special case of \ref{lem:smooth_or_proper:bc} and \ref{lem:smooth_or_proper:hbc}.
    
    The diagonal morphism of $\cs{G}$ is smooth, so that property \bc\ follows from \ref{lem:smooth_or_proper:bc}.
    Note that $\stack G_m(\cs{G}) \cong G^m \times \cs{G}$.
    Thus one sees that the maps $\stack G_m(\cs{G}) \to \stack G_{m-1}(\cs{G})$ are smooth, while the maps $\stack G_m(\cs{G}) \to \stack G_{m+1}(\cs{G})$ are closed immersions (and hence proper).
    Thus property \hbc\ follows from smooth and proper base change.
\end{proof}

\subsection{Some reduction steps}

Let $i\colon\stack Y \to \stack X$ be an open (or closed) immersion.
We have $\ls \stack Y \cong \stack Y \times_{\stack X \times \stack X} \stack Y$ and it follows that there is a canonical open (resp.~closed) immersion of loop stacks $\ls i\colon \ls \stack Y \to \ls \stack X$.
More generally, each $\stack G_m(\stack Y)$ is an open (resp.~closed) substack of $\stack G_m(\stack X)$.

\begin{lemma}\label{lem:lc_immersion}
    Let $i\colon\stack Y \hookrightarrow \stack X$ is a locally closed immersion of QCA stacks.
    If $\stack X$ has property \bc\ or property \hbc, then so does $\stack Y$.
\end{lemma}

\begin{proof}
    We will show the statement for property \bc.
    The proof for \hbc\ is essentially the same.
    First we note that if $\sheaf F \in \catKHol{\stack Y}$, then $i_*\sheaf F \in \catKHol{\stack X}$ since
    \[
      \Hom_{\catDMod{\stack X}}(k_{\stack X}, i_*\sheaf F) \cong
      \Hom_{\catDMod{\stack Y}}(i^*k_{\stack X}, \sheaf F) \cong
      \Hom_{\catDMod{\stack Y}}(k_{\stack Y}, \sheaf F) \ncong 0.
    \]
    Thus it suffices to show that $p_{\stack Y,1,!}p_{\stack Y,2}^!i^! \to \Delta_{\stack Y}^!\Delta_{\stack Y,!}i^!$ is an isomorphism of $\catKHol{\stack X}$.
    We split $i$ into a closed immersion followed by an open one and prove the two cases separately.

    Let us first assume that $i$ is open.
    Using the assumption on $\stack X$ and the fact that $i^! = i^*$ we get
    \begin{multline*}
        p_{\stack Y,1,!} p_{\stack Y,2}^!i^! =
        p_{\stack Y,1,!} (\ls i)^! p_{\stack X,2}^! \cong
        i^! p_{\stack X,1,!} p_{\stack X,2}^! \\ \isoto
        i^! \Delta_{\stack X}^! \Delta_{\stack X,!} =
        \Delta_{\stack Y}^! (i\times i)^! \Delta_{\stack X,!} \cong
        \Delta_{\stack Y}^! \Delta_{\stack X,!} i^!.
    \end{multline*}
    Now let $i$ be a closed immersion.
    As above, it suffices to show that $i_!p_{\stack Y,1,!}p_{\stack Y,2}^!i^! \to i_!\Delta_{\stack Y}^!\Delta_{\stack Y,!}i^!$ is an equivalence.
    Using the assumption on $\stack X$ and the fact that $i_! = i_*$ we get
    \begin{multline*}
        i_! p_{\stack Y,1,!} p_{\stack Y,2}^!i^! =
        p_{\stack X,1,!} (\ls i)_! p_{\stack Y,2}^!i^! \cong
        p_{\stack X,1,!} p_{\stack X,2}^! i_!i^! \\ \isoto
        \Delta_{\stack X}^! \Delta_{\stack X,!} i_!i^! \cong
        \Delta_{\stack Y}^! (i\times i)_! \Delta_{\stack X,!}i^! =
        i_! \Delta_{\stack Y}^! \Delta_{\stack X,!} i^!.
        \qedhere
    \end{multline*}
\end{proof}

Since all operations are local we can check the properties on open covers.

\begin{lemma}\label{lem:open_cover}
    If $\stack X$ has a cover by open substacks $\stack U_i$ with property \bc\ (resp.~with property \hbc), then $\stack X$ has property \bc\ (resp.~property \hbc).
\end{lemma}

\begin{lemma}\label{lem:product}
    If $\stack{X_1}$ and $\stack{X_2}$ both have property \bc\ (resp.~property \hbc), then so does $\stack{X_1} \times \stack{X_2}$.
\end{lemma}

\begin{proof}
    The category $\catDMod{\stack X_1 \times \stack X_2}$ is generated by elements of the form $\sheaf F_1 \boxtimes \sheaf F_2$ with $\sheaf F_i \in \catDModHol{\stack X_i}$.
    For any schematic morphism $f\colon \stack Y \to \stack Z$ of QCA stacks the functor $f^!$ is continuous on $\catDMod{\stack Z}$ and $f_!$ preserves colimits in $\catDModHol{\stack Y}$, since it has a right adjoint.
    Since a full subcategory inclusion reflects all colimits it is thus sufficient to check that the properties holds on elements of the form $\sheaf F_1 \boxtimes \sheaf F_2 \in \catKHol{\stack X_1 \times \stack X_2}$.
    
    Note that $\ls(\stack X_1 \times \stack X_2) \cong \ls \stack X_1 \times \ls \stack X_2$ and more generally $\stack G_m(\stack X_1 \times \stack X_2) \cong \stack G_m(\stack X_1) \times \stack G_m(\stack X_2)$.
    The statement now follows because all functors respect $\boxtimes$.
\end{proof}

We would like to have an extension of Lemma~\ref{lem:product} to fiber products.
Unfortunately, as the discussion in \cite[Section~1.2]{BenZviNadler:arXiv:CharacterTheoryOfAComplexGroup} shows, the category of D-modules on a fiber product $\stack X_1 \times_{\stack Y} \stack X_2$ is rather badly behaved unless $\stack Y$ is a scheme.
Thus we have to restrict ourselves to the following partial result which is sufficient for our applications.

\begin{lemma}\label{lem:product_quotient}
  Let $G$ be an affine algebraic group acting on schemes $X_1$ and $X_2$.
  Suppose the quotient stacks $\stack X_1 = X_1/G$ and $\stack X_2 = X_2/G$ both have property \bc\ (resp.~property \hbc).
  Then so does $(X_1 \times X_2)/G$.
\end{lemma}

\begin{proof}
  The stack $(X_1 \times X_2)/G$ can be rewritten as $\stack X_1 \times_{\cs{G}} \stack X_2$.
  Thus the results of \cite[Section~5.2]{BenZviNadler:arXiv:CharacterTheoryOfAComplexGroup} show that $!$-pullback exhibits the category $\catDMod{\stack X_1} \otimes_{\catDMod{\cs{G}}} \catDMod{\stack X_2}$ as a full subcategory of $\catDMod{\stack X_1 \times_{\cs{G}} \stack X_2}$ containing $\catK{\stack X_1 \times_{\cs{G}} \stack X_2}$.
  This category can be described as the limit of the cosimplicial category with cosimplices
  \[
    \catDMod{\stack X_1 \times \cs{G} \times \dotsb \times \cs{G} \times \stack X_2}.
  \]
  Thus it suffices to show the properties for the stacks $\stack X_1 \times \cs{G} \times \dotsb \times \cs{G} \times \stack X_2$, where they follow from Lemmas~\ref{lem:product} and~\ref{lem:smooth_or_proper}\ref{lem:smooth_or_proper:bg}.
\end{proof}

\subsection{Contractive \texorpdfstring{$\Gm$}{Gm}-actions}

Let $\stack X$ be a QCA stack with an action of $\Gm$.
We call the action \emph{(globally) contractive} if it can be extended to an action of the monoid $\as 1$.
In this case we have in particular an action of the monoid $\{0,1\} \subset \as 1$ on $\stack X$.
Following the convention of \cite[Appendix~C]{DrinfeldGaitsgory:2015:CompactGenerationOfDModOnBunG} we denote by $\stack X^0$ the stack of $\{0,1\}$-equivariant maps $\{0\} \to \stack X$.
The map $\{0\} \to \{0,1\}$ induces a \emph{contraction morphism} $\pi\colon \stack X \to \stack X^0$, while the map $\{0,1\} \to \{0\}$ induces an inclusion morphism $i\colon \stack X^0 \hookrightarrow \stack X$.

We then have the following \emph{contraction principle} (cf.~also \cite[Proposition~3.2.2]{DrinfeldGaitsgory:2014:OnATheoremOfBraden}).

\begin{theorem}[{\cite[Corollary~C.5.4]{DrinfeldGaitsgory:2015:CompactGenerationOfDModOnBunG}}]
    \label{thm:contraction_principle}%
    Let $\stack X$ be a contractive stack with a trivial $\Gm$-action.
    Then there is a canonical isomorphism of functors
    \[
        \pi _! \cong i^! \colon \catDModHol{\stack X} \to  \catDModHol{\stack X^0}
    \]
    given by adjunction
    \[
        \pi _! \to  i^! i_! \pi _!  = i^! (\pi  \circ  i)_! = i^!.
    \]
\end{theorem}

\begin{remark}
    At first glance the requirement that the $\Gm$-action is trivial may seems paradoxical.
    The example that the reader should keep in mind is the stack $\stack X = \as1/\Gm$, where the quotient is taken via the usual action of $\Gm$.
    There is an action of $\Gm$ on $\stack X$ which comes from the action of $\Gm$ on $\as 1$.
    This action is clearly trivial.
    On the other hand it can be extended to an action of $\as 1$, again induced from the action on the cover.
    This action is not trivial and is the one that defines $\stack X^0 = \cs{\Gm}$.
\end{remark}

\begin{remark}
  For our intended application in the proof of Theorem~\ref{thm:torus_quotient}, the special case of this statement for $\Gm$-quotient stacks is sufficient.
  In this case the statement can be proved directly without having to show the full \enquote{stacky} version, cf.~\cite[Proposition~5.3.2]{DrinfeldGaitsgory:2015:CompactGenerationOfDModOnBunG}.
  We have elected to state the theorem as above in order to have a general version of the following proposition for future applications.
\end{remark}

\begin{proposition}\label{prop:contractive}
    Assume that $\stack X$ has a trivial contractive $\Gm$-action such that $i\colon \stack X^0 \to \stack X$ is a closed immersion and both $\stack X^0$ and $\stack U = \stack X\setminus \stack X^0$ have property \bc\ (respectively property \hbc).
    Then so does $\stack X$.
\end{proposition}

\begin{proof}
    Let us first prove the statement for property \bc.
    Let $j\colon \stack U \hookrightarrow \stack X$ be the (open) complement of $i$.
    Consider the distinguished triangle
    \[
        i_!i^! \Delta^!\Delta_! \to \Delta^!\Delta_! \to j_*j^*\Delta^!\Delta_!
    \]
    on $\catDModHol{\stack X}$.
    It suffices to show that the canonical maps on the outside terms to $i_!i^! p_{1,!}p_2^!$ and $j_*j^* p_{1,!}p_2^!$ respectively are isomorphisms.
    
    Using the fact that $j^* = j^!$, base change and the assumption we see that
    \[
        j_*j^*\Delta^!\Delta_! \cong
        j_*\Delta^!\Delta_!j^* \isoto
        j_*p_{\stack U,1,!}p_{\stack U,2}^! j^* \cong
        j_*j^* p_{1,!}p_2^!
    \]
    is an equivalence.

    Write $\pi\colon \stack X \to \stack X^0$ for the contraction map opposite to $i$.
    Then by the contraction principle,
    \begin{equation}\label{eq:lem:bc_contractive}
        i_!i^! \Delta^!\Delta_! =
        i_!\Delta_{\stack X_0}^!(i\times i)^! \Delta_! \cong
        i_!\Delta_{\stack X_0}^!(\pi\times \pi)_! \Delta_! =
        i_!\Delta_{\stack X_0}^!\Delta_{\stack X_0,!} \pi_! \cong
        i_!\Delta_{\stack X_0}^!\Delta_{\stack X_0,!} i^!.
    \end{equation}
    Write $\ls i \colon \ls\stack X^0 \hookrightarrow \ls\stack X$ for the closed inclusion of loop stacks.
    We note that $\ls\stack X$ inherits a contractive $\Gm$ action from $\stack X$.
    Let $\ls \pi\colon \ls\stack X \to \ls\stack X^0$ be the corresponding contraction map.
    Then by assumption and the contraction principle, \eqref{eq:lem:bc_contractive} is equivalent to
    \[
        i_!p_{\stack X_0,1,!}p_{\stack X_0,2}^! i^! =
        i_!p_{\stack X_0,1,!} (\ls i)^! p_{2}^! \cong
        i_!p_{\stack X_0,1,!} (\ls \pi)_! p_{2}^! =
        i_! \pi_! p_{1,!} p_{2}^! =
        i_! i^! p_{1,!} p_{2}^!.
    \]
    Finally, we note that the $\Gm$-action on $\stack X$ induces contractive $\Gm$-actions on each $\stack G_m(\stack X)$ with $(\stack G_m(\stack X))^0 = \stack G_m(\stack X^0)$.
    Thus the statement for property \hbc\ follows in the same way as above.
\end{proof}

\section{Torus quotients}\label{sec:torus}

\begin{theorem}\label{thm:torus_quotient}
    Let $X$ be a normal quasi-projective variety $X$ over an algebraically closed field $k$ of characteristic $0$ with an action of a torus $T$.
    Then the quotient stack $\stack X = X/T$ has properties \bc\ and \hbc.
\end{theorem}

\begin{proof}
  By Sumihiro's theorem \cite{Sumihiro:1974:EquivariantCompletions} and Lemma~\ref{lem:open_cover} we can assume that $\stack X$ is affine with a linear action of $T$.
  By Lemma~\ref{lem:lc_immersion} it further suffices to consider the case that $\stack X = \as{m}/T$ for some $m$.
  By Lemma~\ref{lem:product_quotient}, we can further decompose $\as{m}$ into $T$-eigenspaces.
  Each of those splits into copies of $\as{1}/T$, which in turn can be written as $\as{1}/\Gm \times BT'$ for some subtorus $T'$ of $T$.
  
  By Lemma~\ref{lem:smooth_or_proper}\ref{lem:smooth_or_proper:bg} it remains to consider the stack $\as{1}/\Gm$, where $\Gm$ acts nontrivially.
  If the $\Gm$ action is contractive this case follows immediately from Proposition~\ref{prop:contractive} and Lemma~\ref{lem:smooth_or_proper}.
  Otherwise the inverse of the $\Gm$-action on $\as1$ is contractive and we can use the contraction principle (and hence Proposition~\ref{prop:contractive}) with the monoid $\as[\infty]{1} = \Gm \cup \{\infty\}$ instead of $\Gm \cup \{0\}$.
\end{proof}

Theorem~\ref{thm:main} is now an immediate consequence of Theorem~\ref{thm:torus_quotient} and Corollary~\ref{cor:properties_imply_monads}.

\begin{example}
  To illustrate let us compute $\HCoh(\catDMod{\as1/\Gm})$ where $\Gm$ acts in the usual way.
  Let $i\colon \stack Z \to \ls\stack X$ be the inclusion of the fiber $\stack Z = \Gm \times \cs{\Gm}$ over the origin and $j\colon \stack U \to \ls\stack X$ the inclusion of the open complement $\stack U = \Gm/\Gm = \pt$.
  Let $p\colon \stack X \to \pt$ be the structure map.
  We have to compute $p_*p_{1,!}p_2^!k_{\stack X}$.
  
  First we note that $p_2^!k_{\stack X} = \omega_{\ls \stack X}$.
  Consider now the distinguished triangle
  \[
    p_*p_{1,!}i_!i^!\omega_{\ls \stack X} \to
    p_*p_{1,!}\omega_{\ls \stack X} \to
    p_*p_{1,!}j_*j^!\omega_{\ls \stack X}.
  \]
  The right term is just the cohomology of a point, i.e.~$k$.
  The left term evaluates to
  \[
    \Gamma_c(\Gm, k_{\Gm}) \otimes \GammadR(\cs{\Gm},k_{\cs{\Gm}}) \cong \bigl(k[-1] \oplus k[-2]\bigr) \otimes k[u]
  \]
  with $u$ in degree $2$.
  As the zeroth Hochschild cohomology cannot vanish (it must contain the identity functor), the connecting morphism is zero.
  Hence
  \[
    \HCoh(\catDMod{\as1/\Gm}) \cong k[x]
  \]
  with $x$ in degree $1$.
  More generally
  \[
    \HCoh(\catDMod{\as n/\Gm^n}) \cong k[x_1,\dotsc,x_n].
  \]
  In a future note we plan to give explicit computations for general toric varieties.
\end{example}

\section{Base change and compactification}\label{sec:compactification}

As was noted above, for a general stack there is no reason to expect that the base change morphism is an isomorphism.
In such cases it is then interesting to understand the cone of the morphism, describing how badly base change fails.

\subsection{The general situation}

Let us for the moment return to an arbitrary Cartesian diagram of stacks
\begin{equation}
  \label{eq:c:base-change:square}
    \begin{tikzcd}
        \stack Z \arrow{d}{p} \arrow{r}{q} & \stack X_1 \arrow{d}{f} \\
        \stack X_2 \arrow{r}{g} & \stack Y
    \end{tikzcd}
\end{equation}
with $f$ and $g$ schematic.
As discussed above, in this situation we have a morphism of functors $\catDModHol{\stack X_1} \to  \catDModHol{\stack X_2}$,
\begin{equation}
    \label{eq:c:base-change-morphism}
     p_! q^! \to  g^! f_!
\end{equation}
induced by adjunctions
\begin{equation}
    \label{eq:c:base-change-adjunctions}
    p_! q^! \to 
    p_! q^! f^! f_! =
    p_! p^! g^! f_! \to 
    g^! f_!.
\end{equation}
If $f$ is proper, then \eqref{eq:c:base-change-morphism} is an isomorphism by Proposition~\ref{prop:pre:base-change}.
To understand the behavior for non-proper $f$, we will approximate it by a proper morphism.

\begin{definition}
    A \emph{relative compactification} of a morphism $f\colon \stack X \to  \stack Y$ is a commutative diagram
    \[
        \begin{tikzcd}
            \stack X \arrow[hook]{r}{j} \arrow{dr}[swap]{f} & \bar{\stack X} \arrow{d}{\bar f} \\
            & \stack Y
        \end{tikzcd}
    \]
    where $j$ is an open embedding and $\bar f$ is proper.
\end{definition}

A famous example of such a relative compactification is Drinfeld's compactification of the morphism $\Bun_B \to  \Bun_G$, where $\Bun_G$ is the stack of $G$-bundles on a curve $C$ with $G$ reductive and $B$ is a Borel subgroup of $G$ \cite{BravermanGaitsgory:2002:GeometricEisensteinSeries}.

Let us now assume that in the situation of diagram \eqref{eq:c:base-change:square} there exists a relative compactification of $f\colon \stack X_1 \to  \stack Y$.
Let $\stack X_1^c$ be the closed complement of the open inclusion $j\colon \stack X_1 \hookrightarrow \bar{\stack X}_1$.
Similarly, we denote by $\bar{\stack Z} = \stack X_2 \times _{\stack Y} \bar{\stack X}_1$ and $\stack Z^c = \stack X_2 \times _{\stack Y} \stack X_1^c$ the corresponding fiber products.
The inclusion and projection maps are denoted as is indicated in the following Cartesian diagrams.
\[
    \begin{tikzcd}
        \bar{\stack Z} \arrow{d}{\bar p} \arrow{r}{\bar q} & \bar{\stack X}_1 \arrow{d}{\bar f} \\
        \stack X_2 \arrow{r}{g} & \stack Y
    \end{tikzcd}
    \qquad\qquad
    \begin{tikzcd}
        \stack Z^c \arrow[hook]{r}{i} \arrow{d} & \bar{\stack Z} \arrow{d}{\bar q} \\
        \stack X_1^c \arrow[hook]{r} & \bar{\stack X}_1
    \end{tikzcd}
\]
We note that $\bar{\stack Z}$ is the disjoint union of the closed substack $\stack Z^c$ and the open substack $\stack Z$.
We can now quantify the failure of the base change morphism to be an isomorphism.

\begin{lemma}
    \label{lem:base-change-criterion}%
    In the situation of the Cartesian square \eqref{eq:c:base-change:square}, the cone of the base change morphism $p_! q^! \to  g^! f_!$ is given by
    \[
        \bar p_! i_*i^* \bar{q}^! j_!.
    \]
    In particular, if $i^* \bar{q}^! j_! = 0$, then $p_! q^! \to g^! f_!$ is an isomorphism of functors.
\end{lemma}

\begin{proof}
    Let $\tilde\jmath\colon \stack Z \hookrightarrow \bar{\stack Z}$ be the open inclusion complement to $i$.
    We split the adjunctions in \eqref{eq:c:base-change-adjunctions} in two by using the compositions
    \[
        f = \bar f \circ  j
        ,\quad
        p = \bar p \circ  \tilde\jmath
        \quad\text{and}\quad
        q = \bar q \circ  \tilde\jmath.
    \]
    Thus the adjunction $p_!q^!\to  p_!q^!f^!f_!$ becomes the sequence
    \[
        p_!q^! \to 
        p_!q^! j^! j_! \to 
        p_!q^! j^! \bar f^! \bar f_! j_!.
    \]
    The equality $p_! q^! f^! f_! = p_! p^! g^! f_!$ becomes
    \[
        p_! q^! j^! \bar f^! \bar f_! j_! =
        p_! \tilde\jmath^! \bar q^! \bar f^! \bar f_! j_! =
        p_! \tilde\jmath^! \bar p^! g^! \bar f_! j_!.
    \]
    Finally the adjunction $p_! p^! g^! f_! \to  g^! f_!$ becomes
    \[
        p_! \tilde\jmath^! \bar p^! g^! \bar f_! j_! =
        \bar p_! \tilde\jmath_! \tilde\jmath^! \bar p^! g^! \bar f_! j_! \to 
        \bar p_! \bar p^! g^! \bar f_! j_! \to 
        g^! \bar f_! j_! =
        g^! f_!.
    \]
    Let us apply the same adjunction morphisms in a different order.
    First the recollement adjunctions
    \[
        p_!q^!
        \xrightarrow{\alpha }
        p_!q^! j^! j_!
        =
        \bar p_! \tilde\jmath_! \tilde\jmath^! \bar q^! j_!
        \xrightarrow{\beta }
        \bar p_! \bar q^! j_!,
    \]
    and then the actual base change
    \begin{equation}
        \label{eq:lem:base-change-criterion:pf:split_morphism_base_change}
        \bar p_! \bar q^! j_!
        \to 
        \bar p_! \bar q^! \bar f^! \bar f_! j_!
        =
        \bar p_! \bar p^! g^! \bar f_! j_!
        \to 
        g^! \bar f_! j_!
        =
        g^! f_!.
    \end{equation}
    We note that the adjunction $\alpha \colon \id \to  j^!j_!$ is an isomorphism and the maps in \eqref{eq:lem:base-change-criterion:pf:split_morphism_base_change} compose exactly to the isomorphism of proper base change (cf.~Proposition~\ref{prop:pre:base-change}).
    Thus the cone of the whole composition is the same as the cone of the morphism $\beta $, which is given by the recollement triangle
    \[
        \bar p_! \tilde\jmath_! \tilde\jmath^! \bar q^! j_!
        \xrightarrow{\beta }
        \bar p_! \bar q^! j_!
        \xrightarrow{\phantom{\beta }}
        \bar p_! i_* i^* \bar q^! j_!
        \xrightarrow{+1}.
        \qedhere
    \]
\end{proof}

\subsection{Relative compactifications for quotient stacks}\label{sec:compactification:quotient}%

In the preceding section we simply assumed that a relative compactification of the given map $f$ exists.
We will now construct such a compactification for the diagonal map of a quotient stack.
Thus let $X$ be a scheme of finite type over $k$ and let $G$ be an affine algebraic group over $k$ acting on $X$.
Let $\stack X = X/G$ be the corresponding quotient stack.

Constructing a relative compactification of $\Delta \colon \stack X \to  \stack X \times  \stack X$ is the same as first constructing a $G \times  G$-equivariant relative compactification of $(\proj2, a)\colon G \times  X \to  X \times  X$ (where $a\colon G \times  X \to  X$ is the action map) and then taking the quotient by the $G \times  G$ action\footnote{%
    Here $G \times  G$ acts on $G \times  X$ by $(s_1,s_2) \cdot (g,x) = (s_2gs_1^{-1},\, s_1x)$.
}.
We let
\[
    \Gamma  = \bigl\{(g, x, x, gx) \in  G \times  X \times  X \times  X\bigr\}
\]
be the graph of $(\proj2, a)$.

We pick a $G\times G$-equivariant compactification $\bar G$ of $G$ and let $\bar \Gamma $ be the closure of $\Gamma $ in $\bar G \times  X \times  X \times  X$.
We have an open embedding $j$ of $G \times  X \cong \Gamma $ into $\bar \Gamma $ and a map $f\colon \bar \Gamma  \to  X \times  X$ given by projection on the last two factors.
The composition $f \circ  j$ is equal to $(\proj2, a)$.

Instead of viewing $\Gamma $ as the graph of $(\proj2, a)$ we can drop the third factor and regard $\Gamma $ as the graph of the action map, i.e.
\[
    \Gammasub{a} = \bigl\{(g, x, gx) \in  G \times  X \times  X\bigr\}.
\]
The closure $\barGammasub{a}$ of $\Gammasub a$ in $\bar G \times  X \times  X$ identifies with $\bar \Gamma $.
Thus for ease of notation we will from now on always set $\Gamma  = \Gammasub a$ and $\bar \Gamma  = \barGammasub{a}$.
With this identification, it is clear that $f$ is proper.

\begin{definition}
    Let $\stack X = X/G$.
    With the above construction we set
    \[
        \bar{\stack X} = \rquot{\bar \Gamma }{G\times G}.
    \]
    We have an open embedding $j\colon X \hookrightarrow \bar{\stack X}$ and a proper morphism $\shortbar \Delta \colon \bar{\stack X} \to  \stack X \times  \stack X$ induced by the map $f$ above, such that $\Delta  = \shortbar \Delta  \circ  j$.
\end{definition}

\begin{remark}
    In the case of $G = \Gm$ the compactification $\bar \Gamma $ is explicitly described in \cite{DrinfeldGaitsgory:2014:OnATheoremOfBraden}.
    In particular, if $X$ is smooth it is shown there that $\bar \Gamma $ is smooth over $\bar G = \ps1$.
    It is possible to extend the methods of \cite{DrinfeldGaitsgory:2014:OnATheoremOfBraden} to quotients by higher dimensional tori.
    The resulting constructions are useful for explicit computations.
\end{remark}

\subsection{Base-change of the diagonal}

We will now formulate a criterion analogous to properties \bc\ and \hbc, but using Lemma~\ref{lem:base-change-criterion}.

Let $\stack Y$ be a stack over $\stack X$.
To simplify notation, we set
\[
  \lsY{\stack X} = \stack Y \XXtimes \stack X,
\]
where the fiber product is taken over the maps $\stack Y \to  \stack X \xrightarrow{\Delta } \stack X \times  \stack X$ and $\Delta $.
Further if $f\colon \stack Y_2 \to  \stack Y_1$ is a morphism of stacks over $\stack X$ we write $\ls f$ for the induced morphism $\lsY[\stack Y_2]{\stack X} \to  \lsY[\stack Y_1]{\stack X}$.

We will use Lemma~\ref{lem:base-change-criterion} in order to reformulate the base change criteria in properties \bc\ and \hbc.
To do so we unfortunately have to introduce some additional notation.
Let us fix a compactification $\bar G$ of $G$ and a relative compactification $\stack X \xrightarrow{j} \bar{\stack X} \xrightarrow{\bar \Delta } \stack X \times  \stack X$ as in Section~\ref{sec:compactification:quotient}.
For any $\stack Y \to  \stack X$ this induces a relative compactification $\clsY{\stack X} = \stack Y \XXtimes \bar{\stack X}$.
Let us write $j_{\stack Y}\colon \lsY{\stack X} \hookrightarrow \clsY{\stack X}$ for the open inclusion and $i_{\stack Y}\colon \lscY{\stack X} = \stack Y \XXtimes \stack X^c \to  \clsY{\stack X}$ for the complement.
Thus we have the following two diagrams with Cartesian squares which are central to the remainder of this paper.
The first describes the situation for a loop stack:
\begin{equation}\label{eq:diag:central1}
    \begin{tikzcd}
        \ls\stack X \arrow[hook]{r}{j_{\stack X}} \arrow{d}{p_2} & \cls\stack X \arrow{d}{\bar p_2} & \lsc\stack X \arrow[cm left hook->]{l}[swap]{i_{\stack X}}  \arrow{d} \\
        \stack X \arrow[hook]{r}{j} & \bar{\stack X} & \stack X^c \arrow[cm left hook->]{l}
    \end{tikzcd}
\end{equation}
For the second let $\stack Y_2 \to  \stack Y_1$ be a morphism of stacks over $\stack X$:
\begin{equation}\label{eq:diag:central2}
    \begin{tikzcd}
        \lsY[\stack Y_2]{\stack X} \arrow[hook]{r}{j_{\stack Y_2}} \arrow{d}{\ls f} & \clsY[\stack Y_2]{\stack X} \arrow{d}{{\cls f}} & \lscY[\stack Y_2]{\stack X} \arrow[cm left hook->]{l}[swap]{i_{\stack Y_1}}  \arrow{d} \\
        \lsY[\stack Y_1]{\stack X} \arrow[hook]{r}{j_{\stack Y_1}} & \clsY[\stack Y_1]{\stack X} & \lscY[\stack Y_1]{\stack X} \arrow[cm left hook->]{l}
    \end{tikzcd}
\end{equation}
With this setup we make the following definition.

\begin{definition}\label{def:good}
    We say that a quotient stack $\stack X = X/G$ \emph{\isgood} if the following two conditions hold.
    \begin{enumerate}
        \item $i_{\stack X}^* \bar p_2^! j_!$ vanishes on $\catDModHol{\stack X}$.
        \item For any schematic morphism $f\colon \stack Y_2 \to  \stack Y_1$ of $G$-quotient stacks over $\stack X$ the composition $i_{\stack Y_2}^* (\bar{\ls f})^! j_{\stack Y_1,!}$ vanishes on $\catDModHol{\lsY[\stack Y_1]{\stack X}}$.
    \end{enumerate}
\end{definition}

\begin{proposition}\label{prop:good-is-good}
    If $\stack X = X/G$ \isgood, then there exists a canonical structure of monad on $p_{1,!}p_2^!$ and the morphism $p_{2,!}p_1^! \to  \Delta ^!\Delta _!$ is an isomorphism of monads.
    In particular there is an isomorphism of algebras
    \[
        \HCoh\bigl(\catDMod{\stack X}\bigr)
        \cong
        \opalg{\GammadR\bigl(\stack X,\, p_{1,!} p_2^! k_{\stack X}\bigr)}.
    \]
\end{proposition}

\begin{proof}
  By Lemma~\ref{lem:base-change-criterion}, \goodness\ implies properties \bc\ and \hbc.
\end{proof}

\subsection{Some reductions}

\begin{lemma}
    Let $U$ be a $G$-equivariant open subset of $X$.
    If $X/G$ \isgood\ then $U/G$ \isgood.
\end{lemma}

\begin{proof}
    Set $\stack U = U/G$.
    Let us fix a compactification $\bar G$ of $G$ and let $\bar \Gamma \subset \bar G \times X \times X$ be as before.
    We will only show the second condition of Definition~\ref{def:good}, the first one is similar.
    Thus let $f\colon \stack Y_2 \to  \stack Y_1$ be a map of quotient stacks over $\stack U$ (and hence also over $U$).
    We note that since $\stack U$ is an open subset we have isomorphisms $\lsY[\stack Y_i]{\stack U} \cong \stack Y_i \XXtimes \stack U$.
    Consider the diagram
    \begin{equation}
        \label{eq:lem:cover:diagram}
        \begin{tikzcd}
            \lsY[\stack Y_1]{\stack U} \arrow[hook]{r}{j_{\stack Y_1}^{\stack U}} \arrow[hook]{d}{\alpha } &
            \clsY[\stack Y_1]{\stack U} \arrow[hook]{d}{\beta } &
            \clsY[\stack Y_2]{\stack U} \arrow{l}[swap]{{\cls f}^{\stack U}} \arrow[hook]{d}{\gamma } &
            \lscY[\stack Y_2]{\stack U} \arrow[cm left hook->]{l}[swap]{i^{\stack U}_{\stack Y_2}} \arrow[hook]{d}{\delta } \\
            \lsY[\stack Y_1]{\stack X} \arrow[hook]{r}{j_{\stack Y_1}^{\stack X}} &
            \clsY[\stack Y_1]{\stack X} &
            \clsY[\stack Y_2]{\stack X} \arrow{l}[swap]{{\cls f}^{\stack X}} &
            \lscY[\stack Y_2]{\stack X} \arrow[cm left hook->]{l}[swap]{i^{\stack X}_{\stack Y_2}}
        \end{tikzcd}
    \end{equation}
    The vertical arrows are open embeddings and all squares are Cartesian (where we use the same compactification of $G$ for $\bar{\stack X}$ and $\bar{\stack U}$, the latter being presented by $\bar \Gamma \cap \bar G \times U \times U$).
    Thus
    \begin{equation*}
        i^{\stack U,*}_{\stack Y_2} {\cls f}^{\stack U,!} j_{\stack Y_1,!}^{\stack U} =
        i^{\stack U,*}_{\stack Y_2} {\cls f}^{\stack U,!} j_{\stack Y_1,!}^{\stack U} \alpha ^* \alpha _*=
        \delta ^* i^{\stack X,*}_{\stack Y_2} {\cls f}^{\stack U,!} j_{\stack Y_1,!}^{\stack X} \alpha _*=
        0.
        \qedhere
    \end{equation*}
\end{proof}

A similar argument can be used to reduce the computation to a smooth cover.
In particular we will want to reduce to computations on a cover by schemes.
On such a cover we will use primes to indicate maps on schemes and will use $\sls{}$ instead of $\ls{}$ for the covers of the loop spaces.
If we start with the cover $\bar \Gamma  \to  \bar{\stack X}$ as in the Section~\ref{sec:compactification:quotient}, the pullbacks to the part of diagram \eqref{eq:diag:central1} of interest are
\begin{equation}
    \label{eq:diag:scheme-cover1}
    \begin{tikzcd}
        \Gamma  \arrow[hook]{r}{j'} \arrow{d} &
        \bar \Gamma  \arrow{d} &
        \scls \stack X \arrow{l}[swap]{\bar p_2'} \arrow{d} &
        \slsc \stack X \arrow[cm left hook->]{l}[swap]{i_{\stack X}'} \arrow{d}\\
        \stack X \arrow[hook]{r}{j} &
        \bar{\stack X} &
        \cls \stack X \arrow{l}[swap]{\bar p_2} &
        \lsc \stack X \arrow[cm left hook->]{l}[swap]{i_{\stack X}}
    \end{tikzcd}
\end{equation}
Similarly, the interesting parts of diagram \eqref{eq:diag:central2} become
\begin{equation}
    \label{eq:diag:scheme-cover2}
    \begin{tikzcd}
        \slsY[\stack Y_1]{\stack X} \arrow[hook]{r}{j_{\stack Y_1}'} \arrow{d} &
        \sclsY[\stack Y_1]{\stack X} \arrow{d} &
        \sclsY[\stack Y_2]{\stack X} \arrow{l}[swap]{\scls f} \arrow{d} &
        \slscY[\stack Y_2]{\stack X} \arrow[cm left hook->]{l}[swap]{i_{\stack Y_2}'} \arrow{d}\\
        \lsY[\stack Y_1]{\stack X} \arrow[hook]{r}{j_{\stack Y_1}} &
        \clsY[\stack Y_1]{\stack X} &
        \clsY[\stack Y_2]{\stack X} \arrow{l}[swap]{\cls f} &
        \lscY[\stack Y_2]{\stack X} \arrow[cm left hook->]{l}[swap]{i_{\stack Y_2}}
    \end{tikzcd}
\end{equation}
In both diagrams all vertical morphisms are smooth and the spaces in the top row are schemes.
Explicitly, if $h'\colon Y \to  X$ is a $G$-equivariant morphism of schemes inducing the structure map $\stack Y \to  \stack X$, then
\[
    \sclsY \stack X =
    \biggl\{
        \bigl(g_1,\, y,\, g_2\bigr) \in  G \times  Y \times  \bar G : \bigl(g_2,\, h'(y),\, g_1h'(y)\bigr) \in  \bar \Gamma 
    \biggr\}.
\]

\begin{lemma}\label{lem:base-change:scheme-cover}
    A quotient stack $X/G$ \isgood\ if the following two conditions hold:
    \begin{enumerate}
        \item The composition $(i_{\stack X}')^*(\bar p_2')^! j'_!$ vanishes on $\catDModHol{\Gamma /G\times G}$.
        \item For each equivariant morphism $f'\colon Y_2 \to  Y_1$ of $G$-schemes the composition $(i_{\stack Y_2}')^*(\scls f)^! j'_{\stack Y_1,!}$ vanishes on $\catDModHol{\slsY[\stack Y_1]{\stack X}/G\times G}$.
    \end{enumerate}
\end{lemma}

\begin{proof}
    This follows from the fact that pullback along the smooth vertical morphisms in \eqref{eq:diag:scheme-cover1} and \eqref{eq:diag:scheme-cover2} is conservative \cite[Lemma~5.1.6]{DrinfeldGaitsgory:2013:FinitenessQuestions} and commutes with the other morphisms up to a shift.
\end{proof}

\begin{lemma}
    \label{lem:base-change:cover}%
    If there exists a $G$-stable open cover $U_i$ of $X$ such that all stacks $U_i/G$ \aregood, then $X/G$ \isgood.
\end{lemma}

\begin{proof}
    Let $\stack U_i = U_i/G$ be the corresponding quotient stacks.
    We will again only show the second condition, the first being simpler.

    Let us first show\footnote{
        This is not completely obvious, since the $\bar{\stack U}_i$ do not necessarily form a cover of $\bar{\stack X}$.
        For example, consider $\ps1$ with the usual linear $\Gm$-action and the usual affine cover.
    } that for any $\stack Y \to  \stack X$ the stacks $\lscY{\stack U_i}$ form an open cover of $\lscY{\stack X}$.
    For this it suffices to show that the open subschemes $\sclsY \stack U_i$ cover $\sclsY \stack X$.
    Let $(g_1, y, g_2)$ be a point of $\sclsY \stack X$.
    Then there exists some $U_i$ with $h'(y) \in  U_i$, where $h'\colon Y \to X$ induces the map $\stack Y \to \stack X$ as above.
    But then $g_1h'(y)$ is also in $U_i$ and hence $(h'(y), g_2, g_1h'(y)) \in  \bar{\stack U}_i$.
    Thus $(g_1, y, g_2)$ is in $\sclsY \stack U_i$.

    It now suffices to show that the restrictions of $i_{\stack Y_2}^*\cls{f}^!j_{\stack Y_1,!} \sheaf F$ to $\lscY[\stack Y_2]{\stack U_i}$ vanish for every $\sheaf F \in  \catDModHol{\lsY[\stack Y_1]{\stack X}}$.
    But this follows from the diagram~\eqref{eq:lem:cover:diagram} (for $\stack U_i$ instead of $\stack U$) and the assumption on $\stack U_i$.
\end{proof}

Finally, it can be useful to consider only partial compactifications.
That is, in the constructions of Section~\ref{sec:compactification:quotient}, instead of going all the way to $\bar G$ we only use some subvariety of it.
Thus let $H$ be a $G\times G$-stable subvariety of $\bar G$ containing $G$, and let $^H\bar{\Gamma }$ be the closure of $\Gamma $ in $H \times  X \times  X$.
We set
\[
    ^H\bar{\stack X} = \rquot{^H\bar{\Gamma }}{G\times G}.
\]
Clearly, if $\left\{H_i\right\}$ is an open cover of $\bar G$ by $G\times G$-stable subvarieties (each containing $G$), then $\left\{{}^{H_i}\bar{\stack X}\right\}$ is an open cover of $\bar{\stack X}$.

Let us fix such an open cover.
For any $\stack Y$ over $\stack X$ we obtain corresponding open covers $\left\{ {}^{H_i}\cls_{\stack Y} \stack X\right\}$ and $\left\{{}^{H_i}\lsc_{\stack Y} \stack X\right\}$ of $\clsY{\stack X}$ and $\lscY{\stack X}$ respectively.
We get our usual diagrams for these partial compactifications:
\[
    \begin{tikzcd}
        \stack X \arrow[hook]{r}{{}^{H_i}j} &
        {}^{H_i}\bar{\stack X} &
        {}^{H_i}\clsY[\stack X]{\stack X} \arrow{l}[swap]{{}^{H_i}p_{2}} &
        {}^{H_i}\lscY[\stack X]{\stack X} \arrow[cm left hook->]{l}[swap]{{}^{H_i}i_{\stack X}}
    \end{tikzcd}
\]
and
\[
    \begin{tikzcd}
        \lsY[\stack Y_1]{\stack X} \arrow[hook]{r}{{}^{H_i}j_{\stack Y_1}} &
        {}^{H_i}\clsY[\stack Y_1]{\stack X} &
        {}^{H_i}\clsY[\stack Y_2]{\stack X} \arrow{l}[swap]{{}^{H_i}\cls f} &
        {}^{H_i}\lscY[\stack Y_2]{\stack X}. \arrow[cm left hook->]{l}[swap]{{}^{H_i}i_{\stack Y_2}}
    \end{tikzcd}.
\]

\begin{lemma}\label{lem:base-change:cover-by-relative-compactifications}
    Let $H_i$ be an open cover of $\bar G$ by $G\times G$-stable subvarieties, each containing $G$.
    Then $\stack X = X/G$ \isgood\ if and only if the following two conditions hold for each $i$:
    \begin{enumerate}
        \item The composition $({}^{H_i}i_{\stack X})^*({}^{H_i}\bar p_2)^! ({}^{H_i}j)_!$ vanishes on $\catDModHol{\stack X}$.
        \item For each schematic morphism of $G$-quotient stacks $f\colon \stack Y_2 \to  \stack Y_1$ the composition $({}^{H_i}i_{\stack Y_2})^*({}^{H_i}\cls f)^! ({}^{H_i}j_{\stack Y_1})_!$ vanishes on $\catDModHol{\lsY[\stack Y_1]{\stack X}}$.
    \end{enumerate}
\end{lemma}

\begin{proof}
    Similar to the proof of Lemma~\ref{lem:base-change:cover}.
\end{proof}

\printbibliography
\end{document}